\documentclass{amsart}%
\usepackage{amsmath}
\usepackage{graphicx}%
\usepackage{amsfonts}%
\usepackage{amssymb}

\theoremstyle{plain}

\numberwithin{equation}{section}

\begin{document}
\title{A Note on Applications of Support Vector Machine}
\author{Seung-chan Ahn}
\address{Fermilab, MS 360, P.O. Box 500, Batavia, IL 60510}
\email{ahn@fnal.gov}
\author{Gene Kim}
\address{19 Gage Road, East Brunswick, NJ 08816}
\email{bernhardkim@@yahoo.com}
\author{MyungHo Kim}\thanks{%
The order of authors is alphabetical.}
\address{Genomics Collaborative Inc. 99 Erie Street, Cambridge, MA 02139}
\email{mkim@genecoop.com}


\maketitle

\begin{abstract}
We describe in a rudimentary fashion how \emph{SVM}(support vector machine)
plays the role of classifier in a mathematical setting. We then discuss its
application in the study of multiple \emph{SNP}(single nucleotide
polymorphism) variations. Also presented is a set of preliminary test
results with clinical data.
\end{abstract}

\section{Introduction}

It is a generally accepted wisdom that the causes of biological effects can
be divided into two categories - inhritable(genes from parents) and
environmental(food, gravity, sunlight, surroundings etc). In this paper, we
focus on inheritable factors.  Our suggestion to multiple \emph{SNP}
variations is based on the following general assumptions(For more details,
see \cite{KK}).

\begin{description}
\item[1]  \emph{Suppose all the SNPs are known and there are no
environmental factors. Then each human is determined by a complete set of
SNP variations uniquely.}
\end{description}

Its consequences are: identical twins are exactly the same. Thus it is
possible to classify \emph{SNP} data sets into several subgroups.
Classification(grouping or clustering) is one of basic and important generic
method for distinguishing one from another.

\begin{description}
\item[2]  \emph{To classify objects we are interested in, the most powerful
technique people developed is to numericalize them, in other words, finding
a way of representation into numbers and the collection of numbers into
vectors in a Euclidean space.}
\end{description}

The two assumptions are separately common senses among researchers. The new
twist is that the two assumptions were not considered in the same scope and 
\emph{SVM} offers a powerful machinery to tackle the problem of
classification in a rigorous and systematic way.

\section{Support Vector Machine and Its Analogy}

The concept of \emph{SVM}(Support Vector Machine) was introduced by Vapnik(%
\cite{Va}) in the late 1970's. Since then the idea of \emph{SVM} found its
application in many diverse fields such as machine learning, gene expression
data analysis, high energy physics experiment at \emph{CERN} (European
Organization for Nuclear Research). Why the idea of \emph{SVM} has
been used in such diverse and unrelated fields ? The reason is clear and obvious: \emph{%
SVM}, based on a solid mathematical foundation, attempts to solve a
universal problem of classification, i.e., we need to know which belongs to
which group. The basic idea of \emph{SVM} is deceptively simple. Given a
collection of vectors in $R^n,$ labeled +1 or -1 that are separable by a
hyperplane, \emph{SVM} finds the hyperplane with the maximal margin. More
precisely, the distance between the closest labeled vectors to the
hyperplane is maximal.(Vapnik, cleverly, connected this distance problem to
an optimization problem by using Kuhn-Tucker condition, \cite{Si}). This
hyperplane could be used to determine to which group an unlabeled vector belongs. This machine fits with inductive scientific method.

To give you a definite flavor of \emph{SVM} in everyday experience, let's
consider about familiar concepts , speed limit, height,weight, blood pressure,
Lipid measurements in blood etc. When the speed limit, critical
values for blood pressure of normal people, Lipid mesurements are
determined, people mainly depend on experimenatal data in the past. As a toy
model, we considered an analogy or correspondence between finding the speed
limit on the road and using Support Vector Machine for a criterion to
determine an association between a given set of multiple SNP variations and
a disease or trait.

In mathematical setting, car speed is a point in $R^1,$%
while a set of numbers consisting of SNP variations(or anything we count
several variables at the same time) is represented as a point of $R^n$.

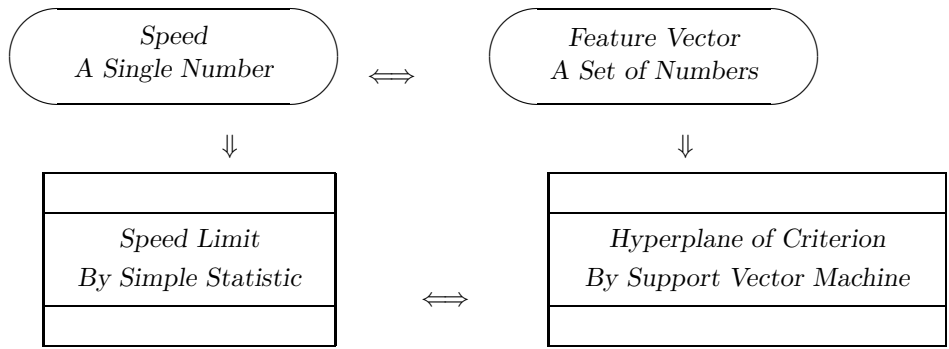
\begin{picture}(124,45)
\put(54,10){\oval(124,36)}
\put(54,10){\makebox(0,15){\sl Speed}}
\put(54,10){\makebox(0,-10){\sl A Single Number}}
\end{picture}
$\Longleftrightarrow $%
\begin{picture}(124,45)
\put(90,10){\oval(124,36)}
\put(90,10){\makebox(0,15){\sl Feature Vector}}
\put(90,10){\makebox(0,-10){\sl A Set of Numbers}}
\end{picture}

\vspace{0.2in}

\hspace{1.0in}$\Downarrow $\hspace{2.3in}$\Downarrow $
\\
\hspace{1.0in}

\begin{picture}(72,45)
\put(5,-15){\framebox(110,50){{\sl By Simple Statistic}}}
\put(5,0){\framebox(110,50){{\sl Speed Limit}}}
\end{picture}
\hspace{1.0in}$\Longleftrightarrow $%
\begin{picture}(72,45)
\put(30,-15){\framebox(150,50){{\sl By Support Vector Machine}}}
\put(30, 0){\framebox(150,50){{\sl Hyperplane of Criterion}}}
\end{picture}

\hspace{1.5in}
\hspace{1.5in}
\\
\\
\begin{description}
\item[Conclusion]  We come to the conclusion that we have to find out a way
of representation of \emph{SNP} variataions at each position. This subject
is open and could be adjusted with experiments for better performance.\footnote{%
After we found out to use the \emph{SVM} to classify multiple SNP
variations, Honki Kim, statistician, pointed out that Classification tree(or
decision tree) might work as well.} Suppose we want to express
east, west, south and north(or DNA letters, \emph{A, C, G, T}). Then we may
represent them as $\{(1,0,0,0),(0,1,0,0),(0,0,1,0)$, $(0,0,0,1)\}$ or $%
\{0.2,0.4,0.6,0.8\}$. This way, at each SNP location, we have a number
depending on genotype in a consistent way, which give us a vector.
\end{description}

\section{Test Results with clinical data}

We generated feature vectors of cardio-patient records by using the same principle described in section 2.
Height, age, sex, weight, ethnic background, medical
history, birth place, blood pressure(systolic and diastolic), Lipid
measurements etc are numericalized and we labeled +1 for a patient who had a history of either heart attack, stroke or heart failure, otherwise -1.  We used Thorsten Joachims' implementation of \emph{SVM}, which gives us the following results(See \cite{Jo}
and, for a different implementation, \cite{Van}). The results strongly
indicate that \emph{SVM} works as intended to separate the data set into
two classes. 

\begin{table}[tbp]
\begin{center}
\begin{tabular}{|l|l|l|l|l|l|l|l|}
\hline
Test & No of Patients & +1 labeled & -1 labeled & C(bound) & Misclassified
& postoneg & negtopos\\ \hline
1 & 1000 & \ \ 212 & \ \ 788 & \ \ 1 & \ \ 56 & \ \ 32 & \ \ 24\\ 
2 & 1000 & \ \ 212 & \ \ 788 & \ \ 2 & \ \ 41 & \ \  23 & \ \ 18 \\ 
3 & 1000 & \ \ 409 & \ \ 591 & \ \ 1 & \ \ 153 & \ \ 37 & \ \ 116 \\ 
4 & 4000 & \ \ 1055 & \ \ 2945 & \ \ 1 & \ \ 438 & \ \ 168 & \ \ 270\\ \hline
\end{tabular}
\end{center}
\caption{Tests with clinical data}
\label{tb:slopes}
\end{table}

For the summary of tests, see the Table 1.

\begin{description}
\item[1] Postoneg means the number of +1 labeled vectors in the group of -1 labeled majority, while negtopos the number of -1 labeled vectors in the group of +1 labeled.

\item[2] Test 1 and 2 are the same data with different C values.
\item[3] Test 1 and 3 are different.
\item[4] Test 3 is contained in Test 4. 
\end{description}
\section{Implication}

Support Vector Machine can be applied for diagnosis of diseases and drug adverse. If, for each possible patient, we input all the test results as a vector, the status of a disease and its prescription could be determined from the past disease records. It should be
noted that the data is not limited to numerical ones and it could include
visual data such as X-ray or MRI image and possibly other sources. For example, in the image data, one extracts area, length, its topological invariant and
others for the totality of input data.

Due to its generic nature of \emph{SVM} already found its application in
diverse field and it may find even more application elsewhere. (Depending on
the users's insight and intuitions, for example, putting genotypes with phenotypes, drug and phenotypes or genotyes etc.)

\section{\textsc{Acknowledgments}}

We are grateful to thank to Profs Larry Shepp of Rugers University, Chul Ahn
of University of Texas at Houston and Dr. Honki Kim of Cedent Corp.
for comments and criticism. Special thanks are due to Prof. Vanderbei at
Princeton University for demonstration of his software implementation of 
\emph{SVM}.

\end{document}